\documentclass[a4paper,12pt]{article}
\usepackage{amsmath,amsthm,amssymb,amsfonts}

\usepackage[latin1]{inputenc}
\usepackage{cite}
\usepackage[english]{babel}
\usepackage[T1]{fontenc}

\newtheorem{teo}{Theorem}
\newtheorem{prop}{Proposition}
\newtheorem{coro}{Corollary}
\newtheorem{lema}{Lemma}

\newcommand{\lrr}{\longrightarrow}
\newcommand{\rr}{\rightarrow}
\newcommand{\inv}[1]{{#1}^{-1}}
\newcommand{\calZ}{{\cal Z}}
\newcommand{\calV}{{\cal V}}
\newcommand{\calT}{{\cal T}}
\newcommand{\jnab}{{\cal J}^\nabla}
\newcommand{\jnabj}{{\cal J}^j}
\newcommand{\hnab}{{\cal H}^\nabla}
\newcommand{\C}{{\mathbb C}}
\newcommand{\Disk}{{\mathbb D}}
\newcommand{\R}{{\mathbb R}}
\newcommand{\dx}{{\mathrm d}}
\newcommand{\zb}{{\overline{z}}}
\newcommand{\XIS}{{\mathfrak X}}

\newcommand{\g}{{\mathfrak g}}       %
\newcommand{\gl}{{\mathfrak{gl}}}     
\newcommand{\m}{{\mathfrak m}}       

\newcommand{\s}{{\mathfrak s}}
\newcommand{\pp}{{\mathfrak p}}
\newcommand{\radicalLie}{{\mathfrak r}}       %
\newcommand{\uni}{{\mathfrak u}}     %

\newcommand{\symp}{{\mathfrak{sp}}}

\newcommand{\na}{{\nabla}}
\newcommand{\Tr}[1]{{\mathrm{Tr}}\,{#1}}

\newcommand{\Ad}[1]{{\mathrm{Ad}}\,{#1}}

\newcommand{\papa}[2]{\frac{\partial#1}{\partial#2}}

\title{\begin{LARGE}\textbf{Remarks on
symplectic twistor spaces}\end{LARGE}}

\author{Rui Albuquerque}

\begin{document}

\maketitle

\vspace{3mm}

\begin{center}
\begin{small}
R. Albuquerque\footnote{The author acknowledges the support of Funda\c{c}\~{a}o para a Ci\^{e}ncia e a Tecnologia, either through POCI/MAT/\-60671/2004 and through CIMA-UE.},
rpa@dmat.uevora.pt, \\
Departamento de Matem\'atica da Universidade de \'Evora \\ 
Centro de Investiga\c c\~ao em Matem\'atica e Aplica\c c\~oes da Universidade de \'Evora (CIMA-UE),\ Rua Rom\~ao Ramalho, 59, 7000 \'Evora, Portugal.

\end{small}
\end{center}

\vspace*{4mm}

\begin{abstract}
We consider some classical fibre bundles furnished with almost complex structures of twistor type, deduce their integrability in some cases and study \textit{self-holomorphic} sections of the general twistor space, with which we define a new moduli space of complex structures. We also recall the theory of flag manifolds in order to study the Siegel domain and other domains alike, which are the fibres of various symplectic twistor spaces. We prove they are all Stein. In the context of a Riemann surface, with its canonical symplectic-metric connection and local structure equations, the moduli space is studied again.
\end{abstract}


\vspace*{2mm}


\vspace*{12mm}

These remarks on complex bundles of complex structures have the purpose of recalling the study of the twistor space of a symplectic manifold $(M,\omega)$, as initiated in \cite{Alb1}, by showing a collection of recent results. The presence of a preferred symplectic connection $\nabla$ is assumed, so we are also studying these differential operators.

Symplectic connections are quite difficult to describe (cf. \cite{Biel}). Also intriguing is the relevance of the symplectic analogue of the Penrose twistor space, which requires a symplectic connection where the latter requires a metric connection. Applications in physics of symplectic twistors, in the way of the celebrated Riemannian case, have not appeared so far in the literature.

In the construction of bundles with complex structures one may admit to have a connection with torsion. Recall we may only have a covariant derivative with, simultaneously, $\nabla\omega=0$ and vanishing torsion if, and only if, $\dx\omega=0$. This is somewhat close to our theorem \ref{generKosMalg}, generalizing another of Koszul-Malgrange since it does not require integrability of the complex structure on the base space.

The pseudo-holomorphic sections $j$ of the (general) twistor bundle are considered in section \ref{Phcs}; satisfying a 1st order equation, we believe they might bring more interesting information in relation with twistors. We found it is possible to define a moduli space of classes of those almost-complex structures $j$, modulo the group of covariantly constant diffeomorphisms.

We prove further details on the manifold topology of symplectic twistor spaces whose sections induce the pseudo-K\"ahlerian metrics on $M$. We recur to the theory of complex symmetric spaces in order to see other characterizations of their standard fibres, which are all Stein domains. With a Grassmannian ambient space it is possible to interchange between signatures. In particular, we prove the fact that, from a symplectic point of view, those bundles with their twistor complex structure all satisfy the respective integrability equations envolving curvature and torsion if, and only if, one of them does.

Finally we initiate some computations for the case of Riemann surfaces and relate the moduli space above to Teichm\"uller space. We give coordinate structure equations of the sheaf of holomorphic functions.

We acknowledge very fruitful conversations with Prof. John Rawnsley, which led to some results here.

\section{Complex bundles of
complex structures and a new moduli space}

\subsection{Twistorial constructions}

A twistor space is an almost complex fibre bundle $\pi:(Z,J^Z)\lrr M$ over a real even dimensional manifold $M$. By this definition, following the one introduced in \cite{Obri}, we mean $Z$ admits a smooth structure of compatible charts and trivializations with a given standard fibre, in the sense of \cite{Steen}, and such that the fibres are almost complex submanifolds. Usually, the complex tangent bundle to the fibres $\ker\dx\pi$ is denoted by $\calV$ due to familiar identification with a \textit{vertical} distribution.

Let $FM$ be the frame bundle of $M$, let $2n$ be the dimension of $M$ and let 
\begin{equation}
 \calZ^0=FM/GL(n,\C)=FM\times_{GL(2n,\R)}GL(2n,\R)/GL(n,\C).
\end{equation}
The projection from the frames $p\in FM$, $p:\R^{2n}\rr T_xM,\ x\in M$, is defined by $p\mapsto pJ_0\inv{p}$, with $J_0$ some fixed linear complex structure of $\R^{2n}$. $\calZ^0$ is the bundle of complex linear structures on the tangents of $M$, called the \textit{general} twistor space of $M$ when it aquires a \textit{Penrose} almost complex structure $\jnab$ induced by a linear connection $\nabla$ on $M$.

Now let $Z$ be as previously. Then we have a commutative diagram
\begin{equation}
 \begin{array}{rcl}
Z & \stackrel{f}\lrr& \calZ^0\\
&\searrow \ \ \:\swarrow & \\
  & M&
\end{array}
\end{equation}
where the map $f$ is defined as the complex structure on $T_{\pi(z)}M$ induced by the isomorphism $\dx\pi:T_zZ/\calV\rr T_{\pi(z)}M$, transporting the complex structure from one space to the other, $\forall z\in Z$. This is, $f(z)$ is the naturally induced complex structure of the quotient space conjugated by the quotient map $\dx\pi$.

There are two very distinct cases that the map $f$ may assume.

One is precisely when $Z$ is a \textit{Penrose} twistor space, a complex subbundle of $\calZ^0$, with the Penrose almost complex structure $\jnab$. Then $f$ is just the inclusion map $Z\hookrightarrow\calZ^0$. The two most familiar subcases are $\calZ^g$ and $\calZ^\omega$, respectively, the orientation and metric $g$-compatible, and the symplectic $\omega$-compatible and tamed complex linear structures on the tangent spaces of, respectively, a Riemannian $(M,g)$ or an almost-symplectic manifold $(M,\omega)$.

Here is the general procedure to construct $\jnab$. The vertical distribution $\calV$ has a natural symmetric space complex structure: vectors $A\in\calV_j,\ j\in\calZ^0$, are endomorphisms of $E=\pi^*TM$ which anti-commute with $j$, hence left multiplication $A\mapsto jA$ is clearly a map with square $-1$. The so called horizontal distribution $\hnab$, suplementary to the vertical one, is given by the kernel of
\begin{equation}\label{projectiontoV}
\pi^*\nabla_{\cdot}\Phi,
\end{equation}
where $\Phi$ is the canonical section $\Phi_j=j\in \mathrm{End}E$. Up to isomorphism $\dx\pi_j:\hnab\rr T_{\pi(j)}M$, each $j$ becomes the \textit{horizontal} complex structure itself, thus living with this ubiquitous character. Adding up, we have $\jnab$ on $T\calZ^0=\hnab+\calV$ preserving such splitting (cf. \cite{BaiWood,BeraOchi,Obri}).

The integrability equations of $\jnab$ are well known: the curvature tensor must satisfy $j^-R(j^+\ ,j^+\ )j^+\:=0,\ \forall j\in\calZ^0$, where $j^\pm=\frac{1}{2}(1\mp ij)$ are the projections inside $TM^c$ onto the eigenspaces of $j$. If $\nabla$ has torsion, then the condition $j^-T(j^+\ ,
j^+\ )=0,\ \forall j\in\calZ^0$, must also be fullfield. Note $R,T$ are real tensors, so if we conjugate the equations we get equivalent statements with $R,T$ again.

The second case we would like to describe arises from a space $Z=Q\times_UN$, where $Q$ is a principal $U$-bundle over $M$ with a $U$-connection, $U$ is a Lie subgroup of $GL(m,\C)$ and $N$ is a $U$-invariant complex manifold. We suppose furthermore that there exists a fibre preserving smooth map $f:Z\rr\calZ^0$. Then $Z$ inherits an almost-complex structure exactly by the same procedure as above: we just keep the vertical complex structure invariantly from $N$ and put $\inv{\dx\pi_|}f(z)\dx\pi_|$ on every horizontal $\hnab_z$, $\forall z\in Z$. Deducing the integrability of this well defined structure seems rather difficult.

The most common situation where the setting occurs is when $M$ itself admits an almost complex structure $j\in\Gamma(M;\calZ^0)$ (denoted by the same letter used for the points in twistor space). The map $f$ in this case is just a constant on each fibre, ie. $f=j\pi$. Moreover, if we want $\pi$ to be a pseudo-holomorphic projection, then $f$ must be $j\pi$. The integrability equations are given next.
\begin{teo}\label{generKosMalg}
Suppose $Z=Q\times_UN$, $f=j\pi$ as above and $Q\subset FM^\C$, the bundle of complex frames. Then $\jnabj$ is integrable if and only if $j$ is integrable and the curvature satisfies $R(j^-\ ,j^-\ )=0$, ie. $R^{(0,2)}=0$.
\end{teo}
\begin{proof}
Notice we just have to prove the result for $Z=Q$ as a principal $GL(m,\C)$-bundle, since the complex connection extends trivially to $FM^\C$, and thence the result follows easily with any $U$-holomorphic invariant factor $N$ (note $U$ may then be a real Lie group). If $Q$ is odd real-dimensional, then one vertical direction is lost when we pass to a quotient, without affecting the induced complex structure.

Let $\theta$ be the soldering 1-form on $Q$, $\theta_p=\inv{p}\dx\pi,\ \forall p\in Q$, and let $\alpha,\rho$ be, respectively, the connection and curvature 1- and 2-forms on $Q$. Defining $\jnabj$ over $Q$ by the procedure described, makes the components of $\alpha,\theta$ become generators of the space of $(1,0)$-forms. Indeed, $\alpha$ is already $\C$-valued and its kernel $\hnab$ is closed under $\jnabj$. For $\theta$ we have
\begin{equation*}
\theta _p\jnabj=\inv{p}\dx\pi\jnabj=\inv{p}j_{\pi(p)}\dx\pi=\sqrt{-1}\theta_p.
\end{equation*}
Now recall the complex structure is integrable iff the space of germs of $(1,0)$-forms generates a $\dx$-closed ideal in $\Omega^*(\C)$. On the principal bundle of \textit{all} frames we have, $\dx\theta=\tau-\alpha\wedge\theta$. With complex frames, we use the projection $j^+:TM^c\rr TM^+$ and find $\dx\theta$ is type $(2,0)+(1,1)$ iff $\tau^{(0,2)}=0$. On the base $M$ this is clearly $j^+T(j^-\ ,j^-\ )=0$. However, this is the Nijenhuis tensor $j^+[j^-\ ,j^-\ ]$. Hence the condition is the same as the integrability of $j$. The well known formula $\rho=\dx\alpha+\alpha\wedge\alpha$ assures us that $\dx\alpha$ has no $(0,2)$ component iff the same happens with $R$.
\end{proof}

Exactly by the same proof but with complex charts in the place of $\theta$, we get the Koszul-Malgrange theorem, which states that any complex vector bundle over a complex ma\-ni\-fold with a complex connection such that $R^{(0,2)}=0$, is a holomorphic vector bundle. Moreover, the latter is the unique holomorphic bundle structure with a same $\overline{\partial}$, namely $\nabla^{(0,1)}$.

This distinct case is indeed less connected with the Penrose twistor construction case than one would think on a first reading. Notice we required a complex connection and complex Lie group rather than the real setting of twistors. One may consult \cite{Obri} to notice the differences in detail.

Finally there is a third case to be considered. Again we assume the manifold $M$ is furnished with an almost complex structure $j$ and a real linear connection $\nabla$ reducible to a principal $G$-bundle of frames $Q\subset FM$. For the purpose, the Lie group $G$ may be understood as a Lie subgroup of the general linear group. We let $J_0\in GL(2n,\R)/GL(2n,J_0)$ be fixed and let $H=\{g\in G:\:gJ_0=J_0g\}$.

Then a twistor space $Z=Q\times_G\frac{G}{H}$ may be defined with an almost complex structure denoted $\jnabj$: the vertical part is given by the usual fibre structure and the horizontal part assumes the fixed $j$, via $\dx\pi_|:\hnab\rr\pi^*TM$. For $G\neq H$, this structure is never integrable.

\subsection{Pseudo-holomorphic complex structures}
\label{Phcs}

We start by a definition. Let $j,j_1\in\Gamma(M;\calZ^0)$ be two almost complex structures on the smooth manifold $M$. A linear connection is said to satisfy condition $A(j,j_1)$ if
\begin{equation}\label{condA}
\nabla_{\XIS^{+,j}}{\XIS^{+,j_1}}\subset{\XIS^{+,j_1}}
\end{equation}
where $\XIS^{+,j_i}$ is the space of $(1,0)$-type vector fields for the structure $j_i$. This definition generalizes the case $A(j,j)$ first studied in \cite{Raw1,Sal2}.
\begin{prop}\label{jholom}
Given a manifold $M$ with a linear connection $\na$ and three sections $j,j_1,j_2:M\rr\calZ^0$, we have:\\
(a) $j_1$ is $(j_2,\jnabj)$-holomorphic if and only if $j_2=j$ and $A(j,j_1)$ is satisfied;\\
(b) in particular, $j$ is $(j,\jnabj)$-holomorphic if and only if $A(j,j)$ is satisfied;\\
(c) $j_1$ is $(j,\jnab)$-holomorphic if and only if $j_1=j$ and $A(j,j)$ is satisfied.
\end{prop}
\begin{proof}
A map $j_1$ is $(j_2,\jnabj)$-holomorphic iff $\dx j_1(j_2X)=\jnabj\dx j_1(X),\ \forall X\in TM$. For the horizontal part we push forward this equation using $\dx\pi$ and get $j_2X=jX$. For the vertical part we apply the projection $\pi^*\nabla_\cdot\Phi$ just to find
\begin{equation*}
 \pi^*\nabla_{\dx j_1(j_2X)}\Phi= \Phi_{j_1}\bigl(\pi^*\nabla_{\dx j_1(X)}\Phi\bigr).
\end{equation*}
This is equivalent to $(j_1^*\pi^*\nabla)_{j_2X}j_1^*\Phi= j_1((j_1^*\pi^*\nabla)_{X}j_1^*\Phi)$. Since $j_1^*\Phi=j_1$ and $j_1^*\pi^*=1$, we get $\nabla_{j_2X}j_1=j_1\,\nabla_{X}j_1$ -- which is equivalent to condition $A(j_2,j_1)=A(j,j_1)$. The rest of the proof is as easy as the previous.
\end{proof}
If $\nabla$ is torsion free, condition $A(j,j)$ implies $j$ integrable. Conversely, if the setting is that of a pseudo-Riemannian triple $(g,\nabla^{LC},j)$ and $\dim M>2$, then the integrability of $j$ implies $A(j,j)$, cf. \cite{Alb2,Sal2}.

Notice $\jnabj$ is never integrable, so we shall follow on with two results on the general well known twistor space. Recall the action of Diff\,$M$ on the space of linear connections, by affine transformation, cf. \cite{Hel1}. A diffeomorphism $g$ of $M$ acts on $\XIS$, the Lie algebra of vector fields, and hence on $\na$:
\begin{equation*}
(g\cdot\nabla)_XY=g\cdot(\nabla_{\inv{g}\cdot X}\inv{g}\cdot Y)\qquad\mbox{where}\quad (g\cdot X)_x=\dx g(X_{\inv{g}(x)}),
\end{equation*}
$\forall X,Y\in\XIS$. The subspaces of torsion free or flat connections are preserved. (On a pseudo-Riemannian setting, by uniqueness, the action of isometries in Levi-Civita is trivial; but in a wider context one must treat it with more circumspect.) The diffeomorphism $g$ also transforms almost complex structures: $j\mapsto \tilde{j}=\dx g\,j\,\dx\inv{g}$. 

In virtue of case c of the proposition above, let us call \textit{self-holomorphic} to those $j$'s which satisfy condition $A(j,j)$.
\begin{prop}
Suppose $j$ is self-holomorphic and $g\cdot\nabla=\nabla$ for some diffeomorphism $g$ of $M$. Then $\tilde{j}$ is self-holomorphic.
\end{prop}
\begin{proof}
It is known that a diffeomorphism $g$ on the base manifold induces a pseudo-biholomorphism $G$ along $g$ on the twistor space. Indeed, $G(j_0)=\dx g\,j_0\,\dx\inv{g}$, $\forall j_0\in\calZ^0$, so $\pi G=g\pi$, and it was proved that $\dx G\,\jnab ={\cal J}^{g\cdot\nabla}\,\dx G$ (cf. \cite{Alb1,Alb2} for two distinct proofs). Now notice $\tilde{j}=G\circ j\circ\inv{g}$. Then 
\begin{eqnarray*}
\dx\tilde{j}_y\,\tilde{j}_y&=&\dx G\,\dx j\,\dx\inv{g}(\dx g\,j_{\inv{g}(y)}\,\dx\inv{g}_y)\\ 
&=& \dx G\,\jnab\dx j_{\inv{g}(y)}\,\dx\inv{g}_y\\
&=& {\cal J}^{g\cdot\nabla}\,\dx G\,\dx j\,\dx\inv{g}_y\ =\ \jnab\, \dx\tilde{j}_y
\end{eqnarray*}
as we wished
\end{proof}
If we work on the symplectic category, then we may prove the converse: if $g$ is a symplectomorphism and $\tilde{j}$ is self-holomorphic, then $g\cdot\nabla=\nabla$. This follows as easily as above, proving $G$ is biholomorphic for $\jnab$ and then applying \cite[corollary 4.1]{Alb1}.

What we have just proved is that there is a good moduli space of self-holomorphic almost complex structures:
\begin{equation}\label{moduli1}
 {\cal M}^\nabla=\frac{\{ j\in\Gamma(M;\calZ^0):\ j\ \mbox{is self-holomorphic}\}}{\mathrm{Diff}(M,\nabla)}
\end{equation}
where Diff$(M,\nabla):=\{g\in \mathrm{Diff}(M):\ g\cdot\nabla=\nabla\}$, the isotropy subgroup of $\nabla$.

For example, let us take a K\"ahler manifold $(M,j_0,\langle\ ,\ \rangle)$ and let $\nabla$ denote the Levi-Civita connection. Then $j_0$ itself shows ${\cal M}^\nabla$ is non-empty.
\vspace*{2mm}\\
\textsc{Remark:} The anti-self-holomorphic almost complex structures $j\in\Gamma(M;\calZ^g)$ in a Riemannian twistor space, endowed with the Levi-Civita connection, correspond to the almost-K\"ahler complex structures in the classification of Gray-Hervella. Indeed, the equation is
\begin{equation}
 j\na_XjY-j\na_{jX}Y+\na_{jX}jY+\na_XY=0
\end{equation}
$\forall X,Y\in TM$. These sections give a moduli space just as above and we recall Diff$(M,\nabla)$ is bigger than the isometry group in general (cf. \cite{Nomi}). Also, notice the twistor space inherits a metric according to its induced tangent bundle decomposition, which, over the image of the embbeding $j$, is given by
\begin{equation}
 g(X,Y)+t\,\Tr{(\na_Xj\,\circ\,\na_Yj)},
\end{equation}
$t\in\R$. The induced K\"ahler form of the second parcel has been studied in the context of K\"ahler-Einstein manifolds (\cite{Seki}) and our methods may relate to the Goldberg conjecture --- any compact almost-K\"ahler Einstein manifold is K\"ahler-Einstein.

\section{On the twistor space of a symplectic manifold}
\label{OttsoaRs}

\subsection{Linear complex structures}
\label{Lcs}

Let us take a path towards the description of the complex structures compatible with a non-degenerate 2-form. Let $V$ be a fixed real vector space of even dimension $2n$ and consider the space of all complex structures in $V$. The action of the real $GL(V)$ shows ${\cal J}(V)=\frac{GL(V)}{GL(V,J)}$.
As a homogeneous space, the tangent space to ${\cal J}(V)$ at $J$ is identified with
\begin{equation}
   \m_{_J}=\left\{ A\in\gl(V):\ AJ=-JA   \right\},
\end{equation}
which is closed under left multiplication by $J$ itself, hence one may prove the space is a complex symmetric space. According to 
\begin{equation}\label{blabla3}
   \gl(V)^c=\gl(V)\otimes\C=\gl(V,J)^c+\m_{_J}^++\m_{_J}^-
\end{equation}
each $A$ equals $(J^+AJ^++J^-AJ^-)+J^+AJ^-+J^-AJ^+$. In sum, ${\cal J}(V)$ has a unique $GL(V)$-invariant complex structure whose $(1,0)$-tangent space at $J$ is $\m_{_J}^+$.

We now specialise to a subspace of ${\cal J}(V)$. Suppose $\omega$ is a symplectic form on the real vector space $V$. Let $J(V,\omega,*)=\{ J\in {\cal J}(V):\ \omega=\omega^{1,1}\ \, \mbox{for}\ J\}$.
Then for any vectors $X,Y$,
\begin{equation*}
0=\omega(X-iJX,Y-iJY)=\omega(X,Y)-\omega(JX,JY)-i(\omega(X,JY)+\omega(JX,Y)),
\end{equation*}
so the new imposed condition is the same as $J$ being a symplectic linear transformation of $V$, or what is called `compatible' with $\omega$. Consider the symmetric form $g_J=\omega(\ ,J\ )$. This non-degenerate inner product has {\it even} signature, say $(2n-2l,2l)$ for some $0\leq l \leq n$, since any maximal subspace where it is positive definite is also $J$-invariant. We denote by $J(V,\omega,l)$ the $n+1$ connected components, as we shall see according to the index $l$, of the disjoint union $J(V,\omega,*)$.
\begin{prop}
(i) None of the $J(V,\omega,l)$ are empty.  \\
(ii) The action of $Sp(V,\omega)$ on $J(V,\omega,l)$ is transitive.
\end{prop}
\begin{proof}
(i) It amounts to show that it is possible to find a basis both symplectic and $g_J$-orthogonal. Recall that a basis $\{X_m\}$ is called symplectic if the matrix of $\omega$ is
\begin{equation*}
 J_0=\left[
\begin{array}{cc}
                0 & -1\\ 1& 0
\end{array} \right] ,
\end{equation*}
so that $X_{m+n}=J_0X_m$ for $m\leq n$ and $J_0$ {\it is} a complex structure.

Let $X,Y\in V\backslash\{0\}$ be such that $\omega(X,Y)=\pm1$ (there exists such a pair). Now, with the non-degenerate $\omega$ restricted, $V_1=\bigl\{U:\ \omega(X,U)=\omega(Y,U)=0\bigr\}=\left\{X,Y\right\}^\omega$
is a symplectic vector space. Assuming the result true by induction for $n-1$, we find the new basis on $V_1$ and a complex structure on $V_1$. Rearranging terms together with $X=-JY$ and $Y=JX$ we find the full basis we required, the index remaining a simple combinatorial problem. Also $V_1$ becomes the orthogonal complement of $\{X,Y\}$. We proved the existence of symplectic bases and the existence of compatible complex structures for any $l$.\\
(ii) Now let $J$ be given and let
\begin{equation}\label{ql}
 Q_l=\left[
\begin{array}{cc}
                1_{n-l,l} & 0 \\ 0 & 1_{n-l,l}
\end{array} \right] .
\end{equation}
$1_{n-l,l}$ is the matrix of the inner product $x_1y_1+\ldots+x_{n-l}y_{n-l}-x_{n-l+1}y_{n-l+1}-\ldots-x_ny_n$. So far we have proved there exists $g\in Sp(V,\omega)$ which transforms, say, the first symplectic basis $\{X_m\}$ into a new one both orthogonal and symplectic. Henceforth $g$ satisfies the equations
\begin{equation*}
       g^tJ_0g=J_0\ \ \ \ \ \ \mbox{and}\ \ \ \ \ \ -g^tJ_0Jg=Q_l.
\end{equation*}
(The necessity of the $-$ sign here becomes obvious when we do $l=0,\
J=J_0$.) Thus $g^tJ_0=J_0g^{-1}$ and hence $g^{-1}Jg=J_0Q_l$. It is trivial to see $J_0$ commutes with $Q_l$, so we have a complex structure $J_0Q_l$ in the same $J(V,\omega,l)$ as $J$ and conjugate to $J$.
\end{proof}
Therefore
\begin{equation}	\label{JVomegal}
  J(V,\omega,l)\ =\ \frac{Sp(V,\omega)}{U(n-l,l)},
\end{equation}
a quotient by the pseudo-unitary group. Fixing any compatible $J$ with index $2l$, we have an inner automorphism $g\mapsto -JgJ$ of $GL(V)$ which preserves $Sp(V,\omega)$ and the respective $U(n-l,l)$. Appealing to the theory we observe that the subspaces we have just been describing are symmetric-subspaces of ${\cal J}(V)$ (cf. \cite{Nomi}). Clearly $J\in\uni(n-l,l)$ too, so there is a direct sum
\begin{equation}	\label{spvomega}
  \symp(V,\omega)=\uni(n-l,l) + \m^s_{_J} 
\end{equation}
where $\m^s_{_J}=\symp(V,\omega)\cap\m_{_J}$. One easily checks that $\m^s_{_J}$ is preserved under left multiplication by $J$, thus the $J(V,\omega,l)$ are also complex submanifolds of ${\cal J}(V)$.

Now we recall the Siegel upper half space or Siegel domain
\begin{equation}	\label{suhs}
 {\cal D}_n=\bigl\{ z \in\C^{\frac{1}{2}n(n+1)}:\ z\ \mbox{symmetric},\ \Im z\ \mbox{positive definite}\bigr\} 
\end{equation}
where the elements $z$ are $n\times n$ matrices with complex entries. $G=Sp(2n,\R)$ acts transitively on ${\cal D}_n$ by
\begin{equation}
  (g,z)\longmapsto g\cdot z=(az+b)(cz+d)^{-1},\ \ \ \ \ \ g=\left[ 
      \begin{array}{cc} a & b  \\ c & d \end{array}  \right] \in G     ,
\end{equation}
where $a,b,c,d$ are square matrices. To see that this is well defined and the action is transitive we appeal to \cite{Siegel}. Easy enough is that the stabiliser of $i1$ is the subgroup of those $g$ for which $d=a$ and $c=-b$, that is, the subgroup $U(n)$. Hence we recover $J(V,\omega,0)$.

Now we look at ${\cal D}_n$ as an open complex manifold.
\begin{prop}						\label{siegel}
The map $\phi:{\cal D}_n\rightarrow J(V,\omega,0)$ given by
\begin{equation}
 x+iy\longmapsto \left[ \begin{array}{cc} xy^{-1} & -xy^{-1}x-y  \\ y^{-1}
& -y^{-1}x \end{array} \right]   
\end{equation}
is a $G$-equivariant {\rm anti}-biholomorphism.
\end{prop}
\begin{proof}
Of course $\phi(i1)=J_0$. Suppose $g\in G$ is such that
\[  g\cdot i1=x+iy. \]
Then $\phi(x+iy)$ must be equal to $gJ_0g^{-1}$. Using the well known
relations between the four squares inside $g$, which also give
$g^{-1}$, we can use the equation above to write $gJ_0g^{-1}$ in terms of $x$ and $y$. Or rather one can check directly that the matrix presented is a true element of $J(\R^{2n},\omega,0)$. By construction, $\phi$ is $G$-equivariant.

Since $G$ acts by rational maps in variable $z$, hence by holomorphic transformations of ${\cal D}_n$, we may conclude that multiplication by $i$ in $T{\cal D}_n$ agrees with a $G$-invariant complex structure. It remains to show that this is the same as {\it right} multiplication by $J_0$ in $\m^s_{_{J_0}}$, up to the iso\-mor\-phism
\[ \dx p_{\mid}: \m^s_{_{J_0}}\longrightarrow T_{i1}{\cal D}_n   \]
arising from the projection $p:G\rightarrow {\cal D}_n$, $\,g\mapsto g\cdot i1$. If we denote
\[    E=\left[ \begin{array}{cc} e_1 & e_2  \\ e_3 & e_4 \end{array}
\right]
\ \ \ \ \ \mbox{and}\ \ \ \ \ g(t)=\mathrm{exp}(tE)=\left[ \begin{array}{cc} a_t & b_t  \\ c_t & d_t
 \end{array} \right],    \]
then the following derivative taking place at point 1 makes sense.
\begin{eqnarray*}
   \dx p(E) &= & \frac{\dx}{\dx t}_{\mid t=0}(a_ti+b_t)(c_ti+d_t)^{-1}   \\
  & = & \dot{a}_0i+\dot{b}_0-i(\dot{c}_0i+\dot{d}_0) =  (e_1-e_4)i+e_2+e_3.
\end{eqnarray*}
Now the result
\[  \dx p(J_0E)\, =\, (-e_3-e_2)i+e_1-e_4\, =\, -i\,\dx p(E)  \]
is immediate to check.
\end{proof}
One may give simple counter-examples to prove that the obvious generalisation of the above to any signature does not hold.

Let us digress in search for a compact symmetric space in which to embed all the $J(V,\omega,l)$. We appeal to the theory of flag manifolds and parabolic subgroups, as explained in \cite{Borel,Burs}. 

Let $G^{C}$ be a connected semisimple complex Lie group and $\g^C$ its Lie algebra. Recall a parabolic subgroup $P$ is the normaliser in $G^C$ of a parabolic subalgebra of $\g^C$, this is, a complex Lie subalgebra $\pp$ which contains a maximal solvable subalgebra of $\g^C$, or a Borel subalgebra. Thus,
\begin{eqnarray*}
 P=\bigl\{ g\in G^C: \ \Ad{(g)}\pp\subset\pp\bigr\}.
\end{eqnarray*}
The theory shows that Lie$(P)=\pp$ and that, if $K$ is a compact real form of $G^{C}$, then
\begin{eqnarray*}
 F=\frac{G^{C}}{P}\simeq\frac{K}{K\cap P}
\end{eqnarray*}
because $K$ acts transitively on $G^{C}/P$ ($K$ has closed image on $F$). These compact and complex spaces are the so called flag manifolds. If $G$ is a non-compact real form (thus semisimple) of $G^{C}$ then it may not act transitively on $F$. On the other hand $G$ is everything but solvable. The \textit{open} orbits of this action are called flag domains.

Applying this theory to $Sp(2n,\C)$ with its two canonical real forms we are able to deduce the components of $J(V,\omega,*)$ are disjoint flag domains in the flag manifold $Sp(n)/U(n)$. We shall prove this in the following lines.

To simplify notation let $V=\R^{2n}$ so that $Sp(V^c,\omega)=Sp(2n,\C):=G^C$. 
Notice $G^C$ acts holomorphically on the Grassmannians $Gr(k,2n)$ of complex $k$-planes (but not transitively). Notice that the stabilizer $P_{\Pi}$ of a point $\Pi$ of the Grassmaniann is also a complex subgroup.

We are particularly interested in the case $k=n$. A quick computation shows that this
\begin{equation}
P=\{g\in Sp(2n,\C):\ g\Pi\subset\Pi\}
\end{equation}
is neither the smallest nor the largest subgroup we can achieve with those actions (perhaps the smallest subgroup is the case when $k=[n/2]$). However, our $P$ is still big enough.
\begin{lema}\label{lema1}
The maximal solvable subalgebra $\radicalLie$ in the Lie algebra $\pp$ of $P$ is maximal solvable in $\symp(2n,\C):=\g^C$, ie. $\radicalLie$ is a Borel subalgebra of $\g^C$.
\end{lema}
\begin{proof}
Just from the theory of solvable Lie algebras over algebraically closed fields, one concludes that all solvable $\s$ in $\gl(\C^{2n})$ preserve some $n$-plane (its elements are all representable in triangular form for a same basis). So any maximal solvable $\s\subset\g^C$ will preserve an $n$-plane. By conjugation with some $g\in GL(2n,\C)$ we find $\s^g\subset\pp$ and $\s^g=\radicalLie$ to be maximal in $\g^C$. To see that inclusion, notice the subgroup $P$ coincides with the stabilizer of the $\radicalLie$-stable $n$-plane.
\end{proof}
Thus, $P$ is a parabolic subgroup and has the `right' dimension. For the computation of the dimension we apply some ideas which J. Rawnsley explained to us.

The theory above gives a framework for studying these problems.
\begin{lema}\label{lema2}
Consider the plane $\Pi=\{(x,0):\ x\in\C^n\}$. Then
\begin{equation} \label{Pigualaetc}
P = \biggl\{ \left[\begin{array}{cc}a & ae\\0 & {a^{-1}}^t \end{array}\right]:\   a\in GL(n,\C),\ e=e^t\in\gl(n,\C)\biggr\}.
\end{equation}
Moreover, we may write $P=GL(n,\C)\rtimes\C^{\frac{n(n+1)}{2}}$.
\end{lema}
\begin{proof}
Let $g\in P$. By definition,
\begin{equation*}
g\left[\begin{array}{c}x \\ 0\end{array}\right]= \left[\begin{array}{cc}a & b \\ c & d \end{array}\right]\left[\begin{array}{c}x \\ 0\end{array}\right]=\left[\begin{array}{c}\ast \\ 0\end{array}\right] 
\end{equation*}
implies $c=0$. Also $g\in G^C$, so
\begin{equation*}
  \left[\begin{array}{cc}a^t & 0 \\ b^t & d^t \end{array}\right] \left[\begin{array}{cc}0 & -1 \\ 1 & 0 \end{array}\right] \left[\begin{array}{cc}a & b \\ 0 & d \end{array}\right] =\left[\begin{array}{cc}0 & -a^t \\ d^t & -b^t \end{array}\right] \left[\begin{array}{cc}a & b \\ 0 & d \end{array}\right] =\left[\begin{array}{cc}0 & -1 \\ 1 & 0 \end{array}\right]
\end{equation*}
and we find $d= {a^{-1}}^t$, \,$d^tb=b^td$. Equivalently, $d= {a^{-1}}^t$ and $b=ae$ with $e^t=e$. 
\end{proof}
Now we take a digression on the real forms of $G^C$, some of which we do not really need. First we see $P\cap U(2n-2l,2l)=U(n-l,l)$. The proof is that
\begin{equation*}
\left[\begin{array}{cc}a^t & 0 \\ (ae)^t & a^{-1} \end{array}\right]Q_l\left[\begin{array}{cc}\overline{a} & \overline{ae} \\ 0 & \overline{{a^{-1}}^t} \end{array}\right]=\left[\begin{array}{cc}a^t1_{n-l,l}\overline{a} & a^t1_{n-l,l}\overline{ae} \\ (ae)^t1_{n-l,l}\overline{a} & (ae)^t1_{n-l,l}\overline{ae}+\inv{a}1_{n-l,l}\overline{{\inv{a}}^t} \end{array}\right]=Q_l
\end{equation*}
if and only if $a\in U(n-l,l),\ e=0$. The matrix $Q_l$ is from (\ref{ql}).
Since $G^C/P$ is connected, arguing with dimensions we find
\begin{equation}\label{embedding0}
\frac{G^C\cap U(2n-2l,2l)}{U(n-l,l)}\ \subset\ \frac{G^C}{P}=\frac{Sp(n)}{U(n)} ,
\end{equation}
where $G^C\cap U(2n-2l,2l),\ \,l\neq0,n$, are the non-compact real forms of $G^C$ and $Sp(n)=G^C\cap U(2n)=U(n,{\mathbb{H}})$ is the well known compact real form of $G^C$. The last equality in (\ref{embedding0}) relies on the fact that the orbit of $Sp(n)$ in $G^C/P$ must be open and closed.

On the other hand, orbits of $Sp(2n,\R)$ are only locally closed. The open ones, the flag domains, certainly appear when $H=P_{\Pi}\cap Sp(2n,\R)$ has the lowest possible dimension as $\Pi$ varies. From \ref{embedding0} above we know this has to be $n^2$. Depending on the signature over $\Pi$ of the metrics given by $h(w_1,w_2)=i\omega(w_1,\overline{w}_2)$, leads to the solutions $H_l=U(n-l,l)$. The $n$-planes are spanned by \begin{equation}\label{nplane}
e_1+if_1,\ldots,e_{n-l}+if_{n-l},e_{n-l+1}-if_{n-l+1},\ldots,e_n-if_n
\end{equation}
where $\{e_j,f_j\}_{1\leq j\leq n}$ is a symplectic basis.

Thus we may claim to have constructed a {\it holomorphic} embedding
\begin{equation}\label{embedding1}
 J(V,\omega,l)=\frac{Sp(2n,\R)}{U(n-l,l)}\longrightarrow \frac{G^C}{P}=\frac{Sp(n)}{U(n)} 
\end{equation}
and a commutative diagram follows:
\begin{equation} \label{maluquice}
\begin{split}
\begin{array}{rcl}
\ \ J(V,\omega,*) & \stackrel{\phi}\longrightarrow & Gr(n,V^c)\ \  \\
\hspace{0cm}\searrow & & \nearrow \\         
  & Sp(n)/U(n) &
\end{array}
\end{split}
\end{equation}
where the map on the top is $J\mapsto V''=V^-$. We shall call \textit{real-Lagrangians} the n-dimensional $\omega$-isotropic $\C$-subspaces $W$ of $V^c$ such that $W\cap\overline{W}=0$.

Notice the first $P$ in lemma \ref{lema2} was the stabilizer of a Lagrangian $n$-plane $\Pi$, since $(x,0)J_0(y,0)^t=(x,0)(0,y)^t=0$, but not a real-Lagrangian.
\begin{prop}
The map $\phi$ is a holomorphic embedding and has image the locally closed manifold $\R Lag(n,V^c)$ of real-Lagrangian subspaces.
\end{prop}
\begin{proof}
Since $\overline{V''}=V'$, the map is injective, and by definition $V''=\phi(J)$ is isotropic. Now let $W\in Gr(n,V^c)$. To any $g\in GL(V^c,\C)$ we associate a sequence
\begin{equation*} 
 W\stackrel{g_\mid}\longrightarrow V^c\stackrel{p}\longrightarrow\frac{V^c}{W}
\end{equation*}
with $p$ only depending on $W$. Therefore $Gr(n,V^c)=\frac{GL(V^c,\C)}{\{\underline{g}=0\}}$
where $\underline{g}=p\circ g_\mid$ and hence
\begin{equation*} 
  T_{W}\,Gr(n,V^c)=\frac{\gl(V^c,\C)}{\{\underline{X}=0\}}\simeq Hom\left(W,\frac{V^c}{W}\right)
\end{equation*}
where $\simeq$ stands for $\underline{X}\simeq p\circ X_{\mid_{W}}$. Now for real $g\in GL(V,\R)$
\begin{equation*} 
 \phi(g\cdot J)=\left\{ v:\ gJg^{-1}(v)=-iv\right\}=gV''.
\end{equation*}
Hence $\dx\phi:\m_{_J}\rightarrow T_{{V''}}\,Gr(n,V^c)$ satisfies $\dx\phi(A)=\underline{A}$ and so
\begin{equation*} 
\dx\phi(JA)=-\dx\phi(AJ)=-p\circ AJ_{\mid_{V''}}=i\,p\circ A_{\mid_{V''}}=i\,\dx\phi(A). 
\end{equation*}
Notice we proved the whole embedding of ${\cal J}(V)$ in $Gr(n,V^c)$ is holomorphic.

Now assume $W\in\R Lag(n,V^c)$. Clearly $\omega:V^c\times V^c\rightarrow\C$ is non-degenerate, so the maximal dimension an isotropic subspace can attain is precisely $n$. Indeed, we have a general formula, $\dim W+\dim W^\omega=2n$, where $W^\omega$ is the $\omega$-anihilator of $W$.

With the above one proves that the hemi-symmetric form on $W$ defined by
\begin{equation*}
  h(w_1,w_2) = i\omega(w_1,\overline{w}_2)= -i\omega(\overline{w}_2,w_1)=\overline{h(w_2,w_1)}
\end{equation*}
is non-degenerate for real lagrangian $W$. According to the signature of this pseudo-metric we may then define $J\in J(V,\omega,l)$ by $Jw=-iw,\ \forall w\in W$, and $J\overline{w}=i\overline{w}$, hence such that $\phi(J)=W$. It is trivial to see $J$ is real. For example
\begin{equation*}
\overline{J}\overline{w}=\overline{Jw}=i\overline{w}=J\overline{w}.
\end{equation*}
We have proved $\phi$ is a biholomorphism onto the aforesaid manifold. Notice also
\begin{equation*}
 \R Lag(n,V^c)\ =\ \{ W:\ W\cap \overline{W}=0\}\cap\{ W: W=W^\omega\} .
\end{equation*}
Here, the first set is open in the Grassmannian and the second is closed.
\end{proof}
The above is only part of either the cell or the algebraic structure of the Gras\-sman\-nian. We will not pursue these in this work. As an example, in $V^c=\C^2$ every line (1-plane) is Lagrangian and there is a circle $S^1$ in $P^1(\C)$ of non-real lines. The open hemispheres are the two Siegel domains ${\cal D}_1=J(\R^2,\omega,0)\simeq\Disk$ and $-{\cal D}_1=J(\R^2,\omega,1)$.

Notice also the following result which was not so clear before, due to phenomena like pseudo-convexity.
\begin{coro}
All $J(V,\omega,l)$ are Stein spaces, $l=0,\ldots,n$.
\end{coro}
\begin{proof}
Conjugating by some non-real $Sp(V,\C)$ element yields a biholomorphism to $J(V,\omega,0)$. The Siegel domain is Stein (it is convex) and being Stein is preserved by biholomorphism, hence the result.
\end{proof}

\subsection{Integrability equations of symplectic twistor space}
\label{Ieosts}

Any given $2n$-dimensional real manifold $M$ endowed with a non-degenerate 2-form $\omega$ has associated to it a bundle $\calZ^{\omega,l}$ of linear complex structures which we call the symplectic twistor space. The standard fibre is $J(V,\omega,l)$, where $V$ is the standard symplectic vector space. We have already mentioned $\calZ^{\omega}=\calZ^{\omega,0}$. 

With a linear connection $\nabla$ such that $\nabla\omega=0$, we may define the \textit{Penrose} almost complex structure $\jnab$ on any symplectic twistor space $\calZ^{\omega,l}$, cf. section 1. The integrability equations were recalled in the same section. It is not obvious that they are equivalent for different $l$, ie. independent of the connected components of $\calZ^{\omega,*}$ (assuming $M$ is connected).
\begin{prop}
The almost complex structure $\jnab$ is integrable on some $\calZ^{\omega,l}$ if, and only if, it is integrable on all.
\end{prop}
\begin{proof}
Let $x\in M$ and $V=T_xM$. Consider first the torsion equation:
\begin{equation}\label{numero1}
J^+T(J^-X,J^-Y)=0  \qquad\forall X,Y\in T_xM\ \ \mbox{and}\ \  J\in \calZ^{\omega,l}_x.
\end{equation}
Fix $J_{0_l}$ in this set. Then (\ref{numero1}) is saying that $T$ takes values in the largest $G=Sp(V,\R)$-invariant subspace of torsion-like tensors such that
\begin{equation}\label{numero2}
 J_{0_l}^+T(J_{0_l}^-X,J_{0_l}^-Y)=0.
\end{equation}
Indeed, since $(gJ_{0_l}\inv{g})^+=g(1-iJ_{0_l})\inv{g}=g J_{0_l}^+\inv{g}$, we have that
\begin{eqnarray*}
J_{0_l}^+(g^{-1}\cdot T)(J_{0_l}^-X,J_{0_l}^-Y) &= &J_{0_l}^+g^{-1}T(gJ_{0_l}^-X,gJ_{0_l}^-Y) \\
 &= & g^{-1}(g\cdot J_{0_l})^+T((g\cdot J_{0_l})^-gX,(g\cdot J_{0_l})^-gY)
\end{eqnarray*}
and hence $T$ is in the subspace iff $\inv{g}\cdot T$ is in the subspace, iff $T$ satisfies (\ref{numero1}) for $J=g\cdot J_{0_l}$. These are ideas from \cite{Obri}.

Notice such subspace of the space of torsion tensors is immersed in a $G^C$-space $\cal T$ of \textit{complex} linear tensors defined by the same condition (\ref{numero2}); such mapping is induced by complexification, ie.
\begin{equation*}
 \{\mbox{all torsion tensors}\}=\wedge^2V\otimes V\longrightarrow\wedge^2V^c\otimes_c V^c.
\end{equation*}
We can now pass to another $J(V,\omega,l')$ by acting on $J_{0_l}$ with an element of $G^C\backslash G$, as we saw in (\ref{embedding1}). Since $\cal T$ is $G^C$-invariant, this space is the same for all $l$ and hence the result.

Finally, analogous arguments follow for the curvature condition, this time with the $G$-subspace sitting in
\begin{equation*}
  \wedge^2V^c\otimes_c S^2V^c 
\end{equation*}
since $\symp(V,\R)=S^2V$.
\end{proof}
We remark additionally that the theory of representations says the irreducible complex $Sp(2n,\C)$-spaces are again $\C$-isomorphic to some $\wedge^qV^c\otimes_c S^pV^c$, with $p,q\geq0$.

As we proved, there is no greater advantage for symplectic geometry in considering the whole $\calZ^{\omega,*}$, --- rather than its 0-component.

\section{The twistor space of a Riemann surface}
\label{TtsoaRs}

Let $(M,g_0)$ be a Riemann surface, where $g_0=h\,\dx z\dx\zb$ denotes the metric in a local conformal coordinate $z$ (this is known to exist, a result due to Gauss). So we are fixing a preferred complex structure on $M$. We have
$h=2g_0(\partial_z,\partial_\zb)>0$, where
\begin{equation}\label{notacao1}
 z=x+iy \hspace{1cm}\mbox{and}\hspace{1cm} \partial_z=\frac{1}{2}\bigl(\papa{ }{ x}-i\papa{ }{y}\bigr)
\end{equation}
and thus $\omega=\frac{ih}{2}\,\dx z\wedge\dx\zb$ is the K\"ahler form.

Now the \textit{symplectic} twistor space of $M$ is the space $\calZ^\omega$ described in previous sections, with fibre identified with the Poincar\'e disk, $SL(2)/U(1)=\Disk$; it admits a real smooth chart $(z,w)$. The second variable appears as follows. Any point $j\in\calZ^\omega$ is given by its $(1,0)$-line in $TM\otimes\C$ (notice the conjugation map is in the $\C$ factor, irrespective of the chosen $j$). Then we pick a standard generator of such line:
\begin{equation}\label{notacao2}
v=\partial_z+\overline{w}\partial_\zb.
\end{equation}
With the chosen chart $z$, the parameter $w$ determines $TM^{1,0,j}$ univocally. As well as the metric $g_j=\omega(\ ,j\ )$ with $g_0$ corresponding to a section $w\equiv0$. Now the condition $-i\omega(v,\overline{v})>0$ gives
\begin{equation}\label{notacao3}
h\,\dx z\wedge\dx\zb(\partial_z+ \overline{w}\partial_\zb,\partial_\zb+w\partial_z) 
=h(1-\overline{w}w)>0\ \ \Longleftrightarrow\ \ |w|^2<1,
\end{equation}
hence the Poincar\'e disk. It is easy to see that a conformal change in coordinates $\papa{ }{z_1}=\papa{z}{z_1}\papa{ }{z}$ yields the change in twistor space coordinates $w_1=w\papa{z}{z_1}/\papa{\overline{z}}{\overline{z}_1}$.

Let us now take the Levi-Civita $\nabla$ for the preferred metric $g_0$. $\na$ is also a complex connection, hence $\na\dx z=\alpha\otimes\dx z$, with $\alpha\in\Omega^{1,0}$. We easily deduce
$\alpha=-\inv{h}\partial h=-\inv{h}\papa{h}{z}\dx z$, just from equation $\nabla g_0=0$ and the fact that $\na$ is real. Moreover, $\na\dx\overline{z}=\overline{\alpha}\otimes\dx \overline{z}$. 

Let us denote by $u$ the unique generator of ${\hnab}^{1,0,j}$, $j=(z,w)$, which projects to $v$:
\begin{equation}\label{notacao4}
u=\partial_z+\overline{w}\partial_\zb+p\partial_w+q\partial_{\overline{w}}
\end{equation}
with $p,q$ to determine. Then
\begin{eqnarray*}\label{conta1}
0\ =\ (\pi^*\na_u\Phi)v &=& \pi^*\na_uiv-\Phi\pi^*\na_uv\\
&=& (i-\Phi)\pi^*\na_{\partial_z+\overline{w}\partial_\zb+p\partial_w +q\partial_{\overline{w}}}(\partial_z+\overline{w}\partial_\zb)\\
&=& (i-\Phi)\bigl(\na_{\partial_z+\overline{w}\partial_\zb} \partial_z+q\partial_\zb +\overline{w}\na_{\partial_z+\overline{w}\partial_\zb}\partial_\zb\bigr)\\
&=& (i-\Phi)\bigl(\na_{\partial_z}\partial_z+\overline{w}\na_{\partial_\zb}\partial_z+ q\partial_\zb +\overline{w}\na_{\partial_z}\partial_\zb+\overline{w}^2\na_{\partial_\zb}\partial_\zb\bigr)\\
&=& (i-\Phi)\inv{h}\biggl(\papa{h}{z}\partial_z+\bigl(hq+\overline{w}^2\papa{h}{\zb}\bigr)\partial_\zb\biggr).
\end{eqnarray*}
The result is again a $(1,0,j)$-vector, hence a multiple of $v$:
\begin{equation*}
 \papa{h}{z}\partial_z+\bigl(hq+\overline{w}^2\papa{h}{\zb}\bigr)\partial_\zb=\lambda\partial_z+\lambda \overline{w}\partial_\zb
\end{equation*}
from which we get the value of $q$:
\begin{equation}
 q=\inv{h}\overline{w}\bigl(\papa{h}{z}-\overline{w}\papa{h}{\zb}\bigr).
\end{equation}
Writting another computation for $(\pi^*\na_u\Phi)\overline{v}=0$, we get $p$:
\begin{eqnarray*}
 0&=& (i+\Phi)\pi^*\na_{\partial_z+\overline{w}\partial_\zb+p\partial_w +q\partial_{\overline{w}}}(\partial_\zb+w\partial_z)\\
&=& (i+\Phi)\bigl(\overline{w}\na_{\partial_\zb}\partial_\zb+w\na_{\partial_z}\partial_z+p\partial_z)\\
&=& (i+\Phi)\bigl(\overline{w}\inv{h}\papa{h}{\zb}\partial_\zb+w\inv{h}\papa{h}{z}\partial_z+p\partial_z\bigr).
\end{eqnarray*}
Again we have an equation
\begin{eqnarray*}
\overline{w}\inv{h}\papa{h}{\zb}\partial_\zb+w\inv{h}\papa{h}{z}\partial_z+p\partial_z&=&\lambda\partial_\zb+\lambda w\partial_z\quad\mbox{for some $\lambda$}
\end{eqnarray*}
\begin{equation}\label{peq}
\Longleftrightarrow\qquad p=\inv{h}w\bigl(\overline{w}\papa{h}{\zb}-\papa{h}{z}\bigr).
\end{equation}
Having found $p,q$ it is convenient to check the integrability of $\jnab$. This will confirm that we made the correct choices of orientation (for it is known that different choices on one direction, horizontal or vertical, induce a non-integrable almost complex structure on twistor space). Indeed, the Lie bracket of the two given independent $(1,0)$-vector fields is
\begin{equation}
 [u,\partial_w]=[\partial_z+\overline{w}\partial_\zb+p\partial_w+q\partial_{\overline{w}},\partial_w]=\papa{p}{w}\partial_w,
\end{equation}
again a $(1,0)$-vector field (since $\partial q/\partial w=0$).
\vspace{3mm}\\
\textbf{Example}. Consider the same setting as above, with $z$ as affine coordinate of projective space $(S^2,g_0,\omega)$ --- Fubini-Study is the round metric. Then $\calZ^\omega_{S^2}$ admits the following non-trivial section, for each $k>0$:
\begin{equation}\label{sectionfors2}
w(z)=\frac{|z|^k}{1+|z|^{2k}}=w(z_1)=\frac{\overline{z_1}^2}{z_1^2}w_1(z_1)
\end{equation}
(invariance of $w$ under coordinate change $z_1=1/z$ being a coincidence). Indeed, $|w|<1,\ \forall z,k$. We remark that in this case $h=\frac{1}{1+|z|^2}$ and thus $p=hw(\overline{z}-\overline{w}z)$.
\vspace{3mm}
\begin{prop}
A function $f\in C^\infty_{{\calZ}^\omega}$ is $\jnab$-holomorphic if and only if
\begin{equation}\label{eq:holofunctiononZ}
 \papa{f}{\overline{w}}=0,\qquad h\papa{f}{\zb}+hw\papa{f}{z}+w\Bigl(\papa{h}{\zb}-w\papa{h}{z}\Bigr)\papa{f}{w}=0.
\end{equation}
\end{prop}
\begin{proof}
 These equations correspond to $\partial_{\overline{w}}(f)=0$ and $\overline{u}(f)=0$.
\end{proof}
\begin{coro}
 A section $j$ of the twistor space, in coordinates $j(z)=(z,w(z))$, is self-holomorphic iff
\begin{equation}\label{eq:self-holo}
w\Bigl(w\papa{h}{z}-\papa{h}{\zb}\Bigr)+h\papa{w}{\zb}+hw\papa{w}{z}=0.
\end{equation}
\end{coro}
\begin{proof}
 $j$ is self-holomorphic iff $f\circ j$ is $j$-holomorphic, $\forall f\in{\cal O}_{\calZ^\omega}$. The result follows from the computation of $h\,\dx(f\circ j)(\partial_\zb+w\partial_z)=$
\begin{eqnarray*}
&=& h\biggl(\papa{f}{\overline{z}}+\papa{f}{w}\papa{w}{\overline{z}}+ \papa{f}{\overline{w}}\papa{\overline{w}}{\overline{z}}+ w\Bigl(\papa{f}{z}+\papa{f}{w}\papa{w}{z}+\papa{f}{\overline{w}}\papa{\overline{w}}{z}\Bigr)\biggr)\\
&=&-hw\papa{f}{z}-w\Bigl(\papa{h}{\overline{z}}-w\papa{h}{z}\Bigr)\papa{f}{w}+h\papa{f}{w}\papa{w}{\overline{z}}+hw\papa{f}{z}+hw\papa{f}{w}\papa{w}{z},
\end{eqnarray*}
which must vanish. We then may devide by the generic $\papa{f}{w}$.
\end{proof}
Now we remark a point $j$ may also be described by
\begin{equation}
 j(\partial_z)=i\frac{1+|w|^2}{1-|w|^2}\partial_z+\frac{2i\overline{w}}{1-|w|^2}\partial_\zb,\qquad j(\partial_\zb)=\frac{-2wi}{1-|w|^2}\partial_z-i\frac{1+|w|^2}{1-|w|^2}\partial_\zb.
\end{equation}
Hence, when $j$ is any section, there is an associated metric on $M$ given by $g_j=\omega(\ ,j\ )$, whose first fundamental form is described in the preferred chart by
\begin{equation}
 g_{j1\overline{1}}:=g_j(\partial_z,\partial_\zb)=h\frac{1+|w|^2}{1-|w|^2}=hl, \qquad g_{j11}:=g_j(\partial_z,\partial_z)=\frac{-2h\overline{w}}{1-|w|^2}=hm.
\end{equation}
Notice $l,m$ are well defined functions of $(z,w)$ and $l+wm=1$.

Recall ${\cal M}^\nabla$ defined in (\ref{moduli1}) and applied to the present situation; it projects to a subspace of the Teichm\"uller space $\calT_M$ of $M$: any class represented by a self-holomorphic $j$ represents a class in $\calT_M$, which is the set of all complex structures modulo the group Diff${}^+(M)$ of orientation preserving diffeomorphisms (cf. \cite{Eells}).

Notice ${\mathrm {Diff}}^+(M,\nabla)$ is not a normal subgroup in the whole transformation group. However, we may still think of comparing those spaces to understand how many distinct classes has a point in $\calT_M$ (we refer to \cite{FarKra,Hel1,Nomi,Tromba} for more details). We recall $T_g{\mathrm {Diff}}^+(M)\simeq\XIS_M$ and $T_g{\mathrm{Diff}}^+(M,\nabla)\simeq\{X\in\XIS_M:\ {\cal L}_X\nabla=0\}$. The latter may be finite dimensional, as the example ${\mathrm {Diff}}^+(\R^2,\dx)=GL(2)\ltimes\R^2$ already computed in \cite{Alb1} shows.

Another problem consists in finding how do self-holomorphic $j$ distribute. The only case we can solve is the following. On a germ of a Riemann surface with metric $h$ (as above) a constant, then the Banach tangent space at $w$ to the space of self-holomorphic maps is given by the space of smooth $\C$-valued bounded functions $F$ such that
\begin{equation}
 \papa{F}{\overline{z}}+\papa{(Fw)}{z}=0
\end{equation}
We find this by differentiating (\ref{eq:self-holo}) and referring to \cite{Eells}. This equation is a variant of the Beltrami equation for quasiconformal mappings.

In the case of the germ of an open subset of $\R^2$ with canonical metric, and therefore with $\nabla=\dx$, we get the tangent space to the $w=0$ section as the set of holomorphic functions $F$, which is a set considerably larger than $T{\mathrm {Diff}}^+(\R^2,\dx)$, as we saw previously.

Henceforth there is hope in finding one non-trivial class of self-holomorphic complex strucures $j$.

\bibliography{twistorcrm}

\begin{thebibliography}{10}

\bibitem{Alb1}
R.~Albuquerque and J.~Rawnsley.
\newblock Twistor theory of symplectic manifolds.
\newblock {\em J. Geom. Phys.}, 56:214--246, 2006.

\bibitem{Alb2}
R.~Albuquerque and I.~Salavessa.
\newblock On the twistor space of pseudo-spheres.
\newblock {\em Differential Geometry and its Applications}, 25:207--219, 2007.

\bibitem{BaiWood}
P.~Baird and J.C. Wood.
\newblock {\em Harmonic morphisms between Riemannian manifolds}.
\newblock LMS Monog. no. 29. Oxford Uni. Press, 2003.

\bibitem{BeraOchi}
L.~B\'erard Bergery and T.~Ochiai.
\newblock On some generalizations of the construction of twistor spaces.
\newblock In T.~J. Willmore and eds. N.~Hitchin, editors, {\em Global
  Riemannian Geometry}, pages 52--58. Ellis Horwood, Chichester, 1984.

\bibitem{Biel}
P.~Bieliavsky, M.~Cahen, S.~Gutt, J.~Rawnsley, and L.~Schwachh{\"o}fer.
\newblock Symplectic connections.
\newblock {\em Int. J. Geom. Methods Mod. Phys.}, 3:375--420, 2006.

\bibitem{Borel}
A.~Borel.
\newblock {\em Linear Algebraic Groups}.
\newblock Springer-Verlag, 2nd edition, 1991.

\bibitem{Burs}
F.~Burstall and J.~Rawnsley.
\newblock {\em Twistor theory for Riemannian symmetric spaces}, volume 1424 of
  {\em Lect. Notes in Math.}
\newblock Springer, Berlin, 1990.

\bibitem{Eells}
J.~Eells.
\newblock A setting for global analysis.
\newblock {\em Bull. Amer. Math. Soc.}, 72:571--807, 1966.

\bibitem{FarKra}
H.~Farkas and I.~Kra.
\newblock {\em Riemann Surfaces}.
\newblock Graduate Texts in Mathematics, Vol. 71. Springer Verlag, 2nd edition,
  1992.

\bibitem{Hel1}
S.~Helgason.
\newblock {\em {Differential geometry, Lie groups, and symmetric spaces.}}
\newblock {Pure and Applied Mathematics, 80. New York-San Francisco-London:
  Academic Press}, 1978.

\bibitem{Nomi}
S.~Kobayashi and K.~Nomizu.
\newblock {\em Foundations of Differential Geometry}, volume 1 and 2.
\newblock Wiley-Interscience, 1969.

\bibitem{Obri}
N.~O'Brian and J.~Rawnsley.
\newblock Twistor spaces.
\newblock {\em Ann. Global Anal. Geom.}, (3):29--85, 1985.

\bibitem{Raw1}
J.~Rawnsley.
\newblock {F-structures, F-twistor spaces and harmonic maps.}
\newblock {Geometry Semin. "Luigi Bianchi", Lect. Sc. Norm. Super., Pisa 1984,
  Lect. Notes Math. 1164, 85-159 (1985).}, 1985.

\bibitem{Sal2}
S.~Salamon.
\newblock {Harmonic and holomorphic maps.}
\newblock {Geometry Semin. "Luigi Bianchi", Lect. Sc. Norm. Super., Pisa 1984,
  Lect. Notes Math. 1164, 161-224}, 1985.

\bibitem{Seki}
K.~Sekigawa and L.~Vanhecke.
\newblock Almost hermitian manifolds with vanishing first chern classes or
  chern numbers.
\newblock {\em Rend. Sem. Mat. Univ. Pol. Torino}, 50(2):195--208, 1992.

\bibitem{Siegel}
C.~L. Siegel.
\newblock Symplectic geometry.
\newblock {\em Amer. J. Math.}, 65:1--86, 1943.

\bibitem{Steen}
N.~Steenrod.
\newblock {\em The Topology of Fibre Bundles}.
\newblock Princeton, 1957.

\bibitem{Tromba}
A.~J. Tromba.
\newblock {\em Teichm\"uller theory in Riemannian geometry}.
\newblock Lectures in Mathematics ETH Z\"urich. Birkh\"auser Verlag, Basel,
  1992.

\end{thebibliography}
\bibliographystyle{plain}

\end{document}